\documentclass[a4paper,12pt]{article}
\usepackage{amsfonts}
\begin{document}

\title{Donald Arthur Preece: A life in statistics, mathematics and music}
\author{R. A. Bailey,\\
School of Mathematics and Statistics,\\
         University of St Andrews,\\
         St Andrews,
         Fife KY16 9SS,
         United Kingdom,\\ and 
 School of Mathematical Sciences (emerita),\\
         Queen Mary University of London,\\
         Mile End Road,
         London  E1 4NS,
         United Kingdom}

\date{}

\maketitle

Donald Arthur Preece died on 6 January 2014, aged 74.  He is survived by his brother Robert.

\section{Edinburgh and St Andrews}
Donald was born in Edinburgh on 2 October 1939.  His father, Isaac Arthur Preece, was a lecturer (later the holder of the first chair) in Brewing and Industrial Fermentations at Heriot--Watt College, which was to become Heriot--Watt University.  His mother, Dorothy Maud, n\'ee Banner, played violin (as a non-professional) in the Reid Symphony Orchestra, which was conducted by Donald Tovey.

From 1945 to 1958 Donald attended George Heriot's School in Edinburgh, where he won medals in both Mathematics and French in 1956.  
He walked to school, and never forgot the day when the Clean Air Act came into force: for the first time, he could see the sea on his daily walk.

He also learnt to play the cello, piano and organ. He was  a cellist in the National Youth Orchestra of Great Britain under the baton of Malcolm Arnold during 1957--1958 and became an Associate of the Royal College of Organists in 1958. 

At the end of school he had to make a choice between his two beloved subjects---mathematics and music.  He decided that he would never be quite good enough to make a decent living as a performing musician, and so studied Mathematics at St Andrews University from 1958 to 1962, obtaining medals in both Pure Mathematics and Applied Mathematics.  

Ian Anderson started his degree in Mathematics at St Andrews in 1960.  He stayed in St Regulus Hall, the same residence as Donald, who offered to be his `senior man' (mentor).  Ian still has the Latin receipt that Donald gave him for a pound of raisins on Raisin Monday.

Donald was part of a cello quartet at St Andrews. He was also a co-founder of the university Madrigal Group, which he directed from 1960 to 1962.

Just as the tug between music and mathematics never left him, so too his later life shows that he was pulled between statistics, specifically the design of agricultural experiments, and combinatorics, specifically defining and constructing new finite combinatorial objects. In reality, these two fields are not so far apart: for example, Latin squares occur in both.

\section{Cambridge and Rothamsted}

After graduating from St Andrews, he spent a year at St Catharine's College, Cambridge. He played in the university orchestra and obtained Cambridge University's Diploma in Mathematical Statistics.  One of his lecturers, David Kendall, encouraged him to seek employment in the Statistics Department at Rothamsted Experimental Station in Harpenden. After remarkably little formal procedure, he took up his first job there, as a Scientific Officer, in 1963.

At that time Rothamsted was one of about forty research stations run by the Agricultural Research Council (ARC): it specialized in agricultural field trials.  Colleagues in the department at the time included H.~Desmond Patterson, John Gower and Gavin Ross.  He also got to know George Dyke, who worked in the Field Experiments Section.
He kept in touch with all of these until the end (George Dyke died in 2012, and Desmond Patterson in 2013). 

Frank Yates was the head of the Statistics Department then.  For the rest of his life, Donald regaled colleagues with stories of FY, as he was known.  He passed on FY's maxim of ``one new idea per paper''; reported that FY claimed to read no papers written by anyone other than R.~A.~Fisher or himself; and retold FY's story that the ``obvious'' idea for balanced incomplete-block designs came to him while he was in the bath.

To researchers used to a university environment, it is almost unbelievable that scientists at Rothamsted could not submit a paper to a journal without first giving copies to both the head of department and the director of the whole of Rothamsted, and obtaining written approval from both.  In the 1980s, I had to obtain such approval from John Nelder and Sir Leslie Fowden for papers on group theory or combinatorial intricacy as well as for those more directly connected with agricultural research. According to Donald, FY would never approve of any paper containing matrices: he disliked them, and could not see the need for them.  However, John Gower reports a less extreme view: although FY himself knew nothing about matrices, John and Michael Healy used them extensively without disapproval.

It was typical for Rothamsted statisticians in the 1960s to publish papers in \textit{Biometrics}, \textit{Biometrika} and \textit{Journal of the Royal Statistical Society, Series~B}.  In 1966 Donald published papers 
 \cite{dap66bcs,dap66rssb,dap66bka} in all three of these, as well as one \cite{dap66tech} in \textit{Technometrics}.  In their different ways, all four of these papers were about a topic that would remain central to Donald's research: designs for experiments in rectangular layouts, and the consequent need to explore properties of relationships between factors, such as rows, columns and various treatment factors.

Rothamsted Manor House had been the home of John Bennet Lawes, who founded  Rothamsted Experimental Station in 1843.
It was bought by the Lawes Agricultural Trust in 1934, and later became effectively a hall of residence for Rothamsted staff and their scientific visitors.  Donald and others now set up the Rothamsted Manor Recitals, which continue to this day.  The musicians at the earliest recitals were all Rothamsted scientists, including Donald himself.  To match the ambience of the Great Drawing Room, where the recitals were held, Donald thought that male members of the audience should wear suits, but few obliged him.

\section{University of Kent}

Donald was promoted to Senior Scientific Officer in 1968.  However, things at Rothamsted were changing.   Frank Yates retired in 1968, 
and Desmond Patterson had moved to the ARC's Unit of Statistics (ARCUS) in Edinburgh in 1967. 
Donald was already in touch 
with S.~Clifford Pearce and 
Geoff(rey) Freeman at East Malling Research Station in Kent, which was also run by the ARC and conducted research on fruit trees.  In 1969 Donald became a Lecturer in Statistics in the Mathematical Institute of the University of Kent at Canterbury, which was in easy reach of East Malling.  In 1972 this university awarded him a PhD in Statistics for published works, and also promoted him to Senior Lecturer.  He was interim head of Statistics for six months in 1973.

Donald lectured on a range of statistical topics, including `Design and analysis of experiments'.  One of his earliest undergraduate students, Miriam Lewis, recalls that he enlivened his teaching by including field trips to East Malling and Rothamsted.

Of course, he did not give up music.  He was conductor of the university's Madrigal Society during 1977--1978.

At Canterbury, Donald collaborated with Clifford Pearce, publishing a paper \cite{dap73bka} with him and J.~R.~Kerr in 1973 on designs for three-dimensional experiments.  This is the first of four ways of generalizing  a rectangular layout by including another factor.  He also maintained contacts with Rothamsted, publishing a paper \cite{dap72gow} with John Gower in 1972 in the \textit{Journal of Combinatorial Theory, Series A}.

Perhaps this is when Donald began to move towards combinatorics.  The British Combinatorial Conference (BCC) began in 1969, and is now held every two years.  Attendance at the first three was by invitation only, but the 1973 conference, in Aberystwyth in July, was open to all, and Donald attended.  He gave a talk about his current research on a problem involving four factors with specified relationships.  The only concise way to describe these relationships is to use matrices. This talk struck a chord with Peter Cameron, who had not met Donald before.  On the conference outing, Peter deliberately sat next to Donald on the bus, and asked him for more information.  There were problems of translation: 
Donald talked about factors on a set of experimental units, but Peter thought about the incidence between the levels of each pair of factors, with no concept of an underlying set;
Donald used lower-case letters for matrices whereas Peter thought that everyone used upper-case letters.  However, they established that Donald's problem was related to something that Peter was thinking about in group theory.  After some more work on both sides, this led to a joint paper 
\cite{dap75cam} in \textit{Utilitas Mathematica} in 1975.  
Thirty years later,
Peter published further work \cite{pjc} on this, and said in the corresponding talk at the 2001 BCC that it was only then that he really understood what Donald had been saying on the bus.

\section{Australia}

From September 1974 to August 1975 Donald took sabbatical leave  from the University of Kent to work in Australia, first in the Department of Mathematical Statistics at the University of Sydney, then in the Division of Mathematics and Statistics of CSIRO (the Commonwealth Scientific and Industrial Research Organisation) in Adelaide, South Australia, where R.~A.~Fisher had spent the last five years of his life  and John Nelder had spent 1965--1966.  Here he worked with W.~Bruce Hall, who specialized in cyclic designs, which are made by developing one or more initial arrays using the integers modulo~$n$. Two papers \cite{dap75oz,dap75hall} in the \textit{Australian Journal of Statistics} followed.  

At this time, statisticians in Adelaide were very receptive to the ideas being put forward by Alan James and Graham Wilkinson \cite{james,jw}, using matrix algebra to explain the relationships between factors.  Of course, this chimed with one of Donald's major interests, and he wrote a paper 
\cite{dap76oz} on non-orthogonal Graeco-Latin designs, which he presented at the Fourth Australian Conference on Combinatorial Mathematics, which took place in Adelaide in August 1975.  Earlier, Clifford Pearce and his colleagues had defined designs according to the pairwise relationships among the factors: a pair may be orthogonal to each other; one may be nested in the other; one may be balanced with respect to the other; and so on \cite{hob,scp}. Now Donald made it absolutely plain that if there are three or more factors whose relationships are neither orthogonality nor nesting then the pairwise relationships are not enough to give the properties of the whole design.  Decades later, I appreciated the clarity of this insight \cite{rab,whatis}.

Donald lived in St Mark's college in Adelaide.  He found that the noisy possums disturbed his sleep, so he used wax ear-plugs.  These melted in his ears, causing him considerable problems.

John Gower was also based temporarily at CSIRO in Adelaide.  One weekend, he and his family took Donald for a trip on a houseboat on the Murray River.  They had a little dinghy to go between the shore and the houseboat.  Donald decided to clean this up, but while doing so he threw away the cork bung.  Fountains of water shot up, while John rushed to get his camera.

\section{Edinburgh again}

Back in the UK by the end of 1975, Donald paid his usual Christmas visit to his mother in Edinburgh (his father had died in August 1964).  As usual on these visits, he called in on ARCUS in early January to see Desmond Patterson.  Here he found me, in the first month of my Science Research Council post-doctoral research fellowship to work on problems of restricted randomization, officially under the supervision of Professor David Finney but in practice under Desmond's guidance.  Donald posed us a problem that he had been thinking about in Australia: if last year's experiment on fruit trees used a Latin square, and the same trees are to be used for this year's experiment, but there may be some residual effects of last year's treatments, how should this year's experiment be designed and randomized?

Donald, Desmond and I met for several hours each day to work on this problem.  Colleagues listened to us arguing over coffee in the common room, at lunch in the staff club, or over a beer.  We obtained good results, which were published \cite{dap78rab} in the \textit{Australian Journal of Statistics} in 1978.  Throughout the discussions, I was in awe of the other two: I could do the necessary permutation group theory, but I did not even know that the word
\textit{bias} had a precise technical meaning.

At the end of Donald's visit, I was astonished when he presented us both with copies of a libretto he had written for an operetta called \textit{When shall we three meet again?}  The reference to the Scottish play was not lost on us.
Desmond was correctly cast as basso profundo, and Donald as tenor.  The only mistake was to portray me as soprano altissima: Donald was not to know that I would go on to sing second alto, or even first tenor, with the Edinburgh University Renaissance Singers.  The text was a fairly close take-off of our actual discussions, gently poking fun at all three of us, but most especially at Donald himself.

\section{Rothamsted again}
In 1978 Donald returned to Rothamsted Experimental Station, this time as Principal Scientific Officer and Biometrics Liaison Officer (Field Crops) of the Ministry of Overseas Development, later the Overseas Development Administration (ODA).  He headed a group of three to four statisticians.  His work included visits to the West Indies, South America, Tanzania, Sudan and Syria, and the British Council sponsored a six-week visit to the University of Ife in Nigeria.

Donald's work for ODA did not prevent him interacting with other Rot\-ham\-sted statisticians.  He was particularly good at nurturing younger colleagues without patronizing them.  For example, he had already done some work himself on identification keys: now he encouraged Roger Payne to get involved, and this led to a joint paper \cite{dap80rss} read to the Royal Statistical Society (RSS) in March 1980.

I joined the Statistics Department at Rothamsted at the start of 1981.  I particularly remember two things that Donald taught me. One was what he called ``sniffing over data'': before you analyse the data, cast your eye over it for anomalies, so that you can ask the scientist to explain these while there is still a chance that he might remember the reason.  Why is that  number double all the rest? Why did the last ten plots to be harvested produce such low yields? Does that pattern of final digits indicate that the data were recorded by various different people? See \cite{dap81st}.

The second thing that Donald taught me was how to read publishers' proofs.  He would hold the publisher's version and read it out, in his immaculate clear enunciation, while I listened and matched it to the original version in front of me.  Both of us would use rulers or cards or fore-fingers to keep our attention on one line of text at a time.

One of the good things about having statisticians in a group, as opposed to being located one per other department, is that you can pick each others' brains over relevant topics: a statistical problem in one area of application may already have been solved in another.  Unfortunately, at this time in Rothamsted there was an unspoken rule that statisticians should not talk about statistics at coffee time or tea time. Rumour had it that one young statistician had once displayed ignorance and that the then head of department, John Nelder, had put him down so sharply (probably intending to denigrate the ignorance rather than the person) that everyone decided to avoid this danger in future.  Donald's solution was to introduce what he called ``topic sessions''.  Every so often the coffee break was turned into a teaching session led by him on a single straightforward topic, such as missing data.  These were always followed by handouts consisting of four sides of A4 paper covered with typescript.

One particularly amusing topic was about Bortkewitsch's horse-kick data. Donald showed the original data, classified by Prussian corps.  Deaths from horse-kicks were not constant across corps: those in the Catering corps suffered relatively few.  With his typical insistence on precise language, Donald complained about those who said that the data were about Prussian officers (\textit{Herren}) when it was actually the foot-soldiers 
(\textit{Heeren}) that died, the officers being safely on horseback.  In his usual generous way, he involved colleagues Gavin Ross and Simon Kirby as co-authors of the ensuing paper~\cite{dap88ross}.

Donald and I often talked to each other about research during these years.  He would say ``I've thought of this new concept: should we call it a 
\textit{genus} or a \textit{species}?'' I would reply ``It doesn't matter.  We're mathematicians.  We can call it anything we like so long as we give a clear definition of the meaning.'' Donald disagreed on two counts.  First, he thought that anyone should  be able to see the pattern in a collection of examples, so that there was no need for a formal definition. Secondly, he thought that the name should help the reader.  More than that, it seemed that helpful terminology was crucial to his own ability to think further about the subject.

One day, he solved a long-standing puzzle of his. He managed to arrange a normal pack of playing cards into a $4 \times 13$ rectangle in such a way that each suit occurred once in each column; each rank (1--10, Jack, Queen and King) occurred once in each row; each row had four cards of one suit and three of each of the others; and each pair of ranks occurred together in precisely one column: see \cite{dap82uti}.  He glued the cards in this arrangement onto a cardboard backing, framed this, and hung it up in the department, in the Fisher building.  I suggested that he should call this arrangement, and its generalization to $k \times v$ rectangles, a 
\textit{double Youden rectangle}.  For once, Donald was delighted with my terminology: he claimed that this name enabled him to think clearly about the structure, and went on to write several papers on the subject.  Many years later, when the Fisher building was closed, Gavin Ross rescued the framed design and returned it to Donald.

During his second period at Rothamsted, Donald continued to maintain contact with Clifford Pearce and Geoff Freeman, particularly over designs for rectangular layouts. The latter, who had now moved to the ARC's National Vegetable Research Station at Wellsbourne, in Warwickshire,  wrote two papers \cite{free,free2} about quasi-complete Latin squares, which have the neighbour-balance property that each unordered pair of letters occurs as neighbours twice in rows and twice in columns.  This inspired me to take the work further in \cite{terr}, and this in turn seems to have inspired much of Donald's combinatorial work after 2000.  

In another generalization of rectangular layouts, each plot is divided into subplots.  Semi-Latin squares are suitable designs for such layouts.
Donald  and Geoff Freeman also published a paper \cite{dap83free} on these, which I and my students drew on for later work.

However, during this period Donald's papers were less about new ideas.  He wrote a series of what were essentially literature reviews 
\cite{dap75oz,dap77uti,dap79oz,dap82uti3}.  With his liking for precise terminology, he was upset that different statisticians could use words like \textit{orthogonal}  or \textit{balanced} to mean different things.  
He compiled a bibliography of randomization, and kept it up to date, but never published it. He did publish expository articles 
\cite{dap80st,dap81st,dap82st,dap86st,dap87st}
in the \textit{Statistician} about simple parts of statistics and how to teach them, as well as some rants 
\cite{dap82uti2,dap82hey,dap84bcs,dap86ag,dap87st2} about statistical practice in the design of real experiments and the subsequent data analysis.
These may have been motivated by poor practice that he found in experimental stations on his overseas trips or by the struggles of his junior colleagues.

John Nelder had established the statistical programming language Genstat at Rothamsted.  Donald soon learnt how to use it, and published a series of notes 
\cite{dap84gen,dap86gen,dap87gen,dap88gen,dap94gen,%
dap94gen2,dap95ross,dap97ross} 
in the \textit{Genstat Newsletter} about his favourite designs.

\section{East Malling}

In 1975 Clifford Pearce retired from East Malling and became Professor of Biometry at the University of Kent. He was succeeded by Ken Martin.
In 1985, Donald took his position as Head of the Statistics Department at East Malling Research Station.  Because of the natural rectangular layouts, designing experiments on orchards was very congenial to him.  He also had an excellent junior colleague in Martin Ridout.  He moved back to Kent, and the University of Kent made him an honorary Reader in Agricultural Biometrics in 1987.

At East Malling, Donald became involved with the Luncheon Club. Although this was technically open to all staff, its members were almost entirely retired.  He organized several musical events there, using the grand piano in Bradbourne House. He frequently dressed up in special costumes, and expected other performers to do likewise. Even as a humble page-turner, Martin had to wear fancy dress.

Sadly, things did not work out professionally for Donald at East Malling as well as he had hoped.  Clifford Pearce had earnt enormous respect from his scientific colleagues during his decades there.  Donald had expected to walk into his shoes, and was taken aback to discover that he would have to earn respect himself afresh.  Moreover, the environment had changed in the ten years since Clifford had left.   The emphasis on field experiments was reduced, and Donald had to oversee  the closing of the Records Office, which had previously kept details of all trials ever conducted there.  Given his frequent use of the archives while at Rothamsted, this must have upset him.
There was a series of changes in the ARC, which renamed itself the Agriculture and Food Research Council (AFRC) in 1983 and became part of the new Biotechnology and Biological Sciences Research Council in 1994. The AFRC merged its previously separate research stations into ten institutes, and re-organized personnel, partly under the common misapprehension that there was less need for statisticians now that all scientists had computers on their desks.  In 1990 Donald 
took redundancy from
what was now called the Institute of Horticultural Research at East Malling, as did many colleagues. Martin Ridout stayed until 2000 before joining the University of Kent.

Later, Donald was to say that several other redundant scientists of his own age simply did not know what to do with themselves.  He described it as ``enforced retirement'', and said that he never wanted to be retired again.

\section{University of Kent again}

Fortunately, the University of Kent offered him a position, first as a part-time lecturer, then (from 1994) as a half-time Professor of Statistics.  He threw his energies into teaching again, especially project students.  He encouraged junior colleagues, advising them that when they finished writing a paper they should put it into a drawer for a while rather than submitting it immediately.  He supervised his first and only PhD student, Chris Christofi, who completed his thesis  \cite{thesis}
in 1993, and went on to publish two papers \cite{cc1,cc2} on the topic.  Donald began new collaborations in combinatorics with his colleague Barry Vowden, with local retiree David Rees,  
and with Peter Owens at the University of Surrey, who retired in 1993 but remained active.

Most exciting was an unforeseen new collaboration that developed.  
In \cite{jan}, John Nelder had been concerned with the problem of computing the matrix product $XX^\top$ when only $k$ of the $v$ rows of $X$ could be held in memory at any one time.  He sought a linearly ordered collection of $b$ blocks of $k$ row-numbers such that (i)~each pair occurs in at least one block and that (ii)~each block differs from its predecessor in just one element.  Donald's first overtly combinatorial publication \cite{dap72gow} was a step towards solving this problem, with the extra condition that
(iii)~the new element in each block does not produce any pairs that have occurred before.

Robin Constable was a combinatorial statistician at the University of St Andrews, who had visited East Malling during his statistical training in Aberdeen.  He had met Donald at the Aberystwyth meeting in 1973, and stayed in touch.  He spent the first three months of 1978  on sabbatical leave at the University of Kent. He worked with Donald, finding some errors in 
\cite{dap72gow} and discovering more designs satisfying (i)--(iii), but nothing was published.

In 1992 Donald was an invited speaker at the Carbondale Combinatorial Conference at Southern Illinois University (SIU) and the following meeting of the American Mathematical Society at Dayton, Ohio.  At a party in Wal Wallis's house in Carbondale, he met Guo-Hui (Grant) Zhang, who asked him what he would speak about in Dayton.  Donald told him about the work with Robin Constable, and was astonished to find that Grant and Wal had two papers \cite{wal,zhang} in press on the same subject, having been posed essentially the same problem by an electrical engineer.  They had defined 
tight (for condition~(iii)) single-change (condition~(ii)) 
covering (condition~(i)) designs. This led to an exciting collaboration involving three further people from SIU.

The story is told in \cite{dap95ica}.  Not only was this collaboration exciting in its own right.  Donald went on to write several further papers, on various topics, with his new collaborators.

However, Donald did not break his links with agricultural statistics.  He took part in the Genstat conference held at Canterbury in July 1993, where he and John Nelder played piano duets at the concert that he organized.

In 1996 John Morgan, in the USA, published a joint paper \cite{jpmudd} and a handbook chapter \cite{jpm} about nested designs: blocks are divided into sub-blocks, and the properties of the design in whole blocks and the design in sub-blocks are both important.  Alternatively, 
in the third generalization of rectangular layouts,
the blocks may be rectangles, and the properties of the design in blocks, in rows and in columns are all important.  This reminded Donald of work \cite{dap67bka}  that he had published in \textit{Biometrika} in 1967, as well as a project that he was then undertaking with David Rees, so he contacted John, and two joint papers \cite{dap99jpm,dap01disc} followed. 

Donald became a regular attender at the BCC, where he gave himself the task of organizing the concert and playing the piano for anyone who needed accompaniment.  At short notice, the promised venue for the 1999 BCC was withdrawn.  Donald stepped into the breach, and offered to host the 1999 BCC at the University of Kent, being joint local organizer with his colleague John Lamb.  This offer was accepted, and so Donald became a member of the British Combinatorial Committee for two years.  Unfortunately, John Lamb became seriously ill in early 1999, and Donald had to do the lion's share of running a very successful conference.

\section{Learned societies}

Throughout these working years, Donald was active in the International Biometric Society (IBS), the RSS, and other learned societies. He  
joined the British Region of the IBS.
He was elected a Fellow of the RSS in December 1963,  an Ordinary Member of the International Statistical Institute in 1977, and a Fellow of the Institute of Statisticians in 1982. He joined the Institute  of Combinatorics and its Applications in 1995.

In the RSS, he encouraged colleagues to attend read papers and contribute to the discussion. He
was a member of the General Applications Section Committee, 1970--1971; of the Journals Committee, 1975--1978; of Council, 1976--1980; of the Local Groups Coordinating Committee, 1978--1982, being chairman 1978--1980; of Programme Committee, 1978--1979 and 1993--1996; of the Research Section Committee, 1983--1986; of the Editorial Committee, 1993--1995; and the Editorial Policy Board, 1996.

In the British Region of the IBS he was keen on the sessions for young biometricians, and may even have been instrumental in setting them up.  As a PhD student, Steven Gilmour talked in one of these sessions in the year 1989--1990.  After the talk, Donald told him about Hadamard matrices of 
order~$16$, and posed some questions about their use for non-regular factorial designs.
He attended some of the international conferences (IBCs) of the IBS.  He took his mother to the one in Constanza, Romania in August 1974, where he booked separate rooms in the names of Mr and Mrs Preece.  The hotel was short of rooms, and thought that they ought to share.  At the IBC  in Guaruja, Brazil in August 1979, the ODA insisted that he use a cheaper hotel than the standard one on offer: Donald reported that he found himself in what was essentially a brothel for gay seamen.
He was a member of the Editorial Advisory Committee of the IBS, 1985--1989; and the British Region Committee, 1987--1990.

He served as a member of the Overseas Committee of the Institute of Statisticians, 1982--1984; and as a member of Council of the Institute of Combinatorics and its Applications, 1999--2002.

He also undertook a wide range of editorial responsibility.  He was British Regional editor of \textit{Statistical Theory and Method Abstracts} from 1966 to 1974. He was an associate editor of \textit{Biometrics} from 1979 to 1989, and Abstracts editor for the \textit{Biometric Bulletin} 
from 1988 to 1991.  He was a member of the editorial boards of the 
\textit{Journal of Agricultural Science} from 1988 to 1994, and of 
\textit{Utilitas Mathematica} from 1990 to 2000, then Managing Editor from 2001 to 2005.  

Most importantly, he was joint editor of \textit{Applied Statistics} from 1993 to 1996, sharing the editorship first with Wojtek Krzanowski and then with Sue Lewis.  Wojtek had overlapped with him for a year at Rothamsted in 1968--1969, where he had benefited from the help and respect that Donald always showed to newcomers.
A published letter from the joint editors in 1994 reminded prospective authors that it was not enough for a new statistical method to be applicable in many areas; to be suitable for publication in \textit{Applied Statistics}, a paper must be motivated by a real application in a single area.  Although this letter had been prompted by the large number of submitted papers containing little evidence of any practical application,
it did seem a rather strange contrast to the direction of Donald's own research interests at the time.

When Sue Lewis replaced Wojtek she benefited greatly from Donald's kind
advice and support. He was devoted to the journal, and gave a tremendous
amount of his time to explaining to authors how their papers could be
improved, or why they were better suited to other journals. In spite of this, Byron Morgan reported that this work load did not seem to distress Donald at all.  Donald strongly believed that the Summary at the start of each paper should be a statement of what the authors had found out, not merely of what they were investigating, and he campaigned enthusiastically for the summaries to match his view.

Just as his music had not been abandoned, neither had his French.  He spent January 1973 as exchange professor at the Orsay campus of the Universit\'e de Paris XI (Paris-Sud).  
It used to be that French and English were equally acceptable languages to the IBS.  In September 1982 the IBC was held at Toulouse, in France.  Donald spoke on \textit{La biom\'etrie: pas rites mais science}.  In September 1991 he was one of several people giving half-day sessions at the European Community Advanced Course in Statistics (Design of Experiments) at S\`ete, in France. Speakers had all been firmly told that this was an international meeting and so they must speak in English.  However, Donald observed that most of the participants were francophone, so he conducted his session in French.

Donald was a pianist and organist, and surely practised before giving performances.  His lectures, whether in English or French, also seemed like polished performances. I was never sure whether he had learnt his script by heart or had simply rehearsed giving the lecture.  It is true that informal questions in the middle could sometimes throw him completely off track.

\section{Queen Mary, University of London}

Donald's mother died in September 1997.  That gave him a certain amount of financial independence, as well as personal independence.  He told Byron Morgan, who was then head of the Institute of Mathematics and Statistics at the University of Kent, that he wanted to resign his position on his 
60th birthday, in 1999, to give himself time to do the things that he really wanted to do.  This was granted, and he was made honorary Professor of Combinatorial Mathematics.

I had moved to the University of London in 1991, first at Goldsmiths College, then at Queen Mary and Westfield College, now called Queen Mary, University of London (QMUL).  I had been running weekly seminars on Design of Experiments during that time, and Donald had attended those when his teaching commitments allowed.  Peter Cameron had been running the flourishing Combinatorics Study Group at QMUL since 1986.  I therefore persuaded QMUL to offer Donald a professorial fellowship at QMUL from 2000.  The salary was derisory: the purpose was mutual benefit to him and us from research interaction.

It seems appropriate to quote here part of the case that I made for his appointment. ``One of Preece's strengths is as a collaborator and inspirer.  By asking awkward questions, he goads other people into doing research.
\ldots\  He is also very good at encouraging both PhD students and staff to keep going on a piece of research that had appeared to be stuck, or to embark on new projects.'' 

Donald accepted the position at QMUL.  He continued to live at his house near East Malling, and also maintained and used his mother's old house in Edinburgh. He undertook some teaching at the University of Kent, especially the supervision of final-year undergraduate projects.  One of these led to a four-author paper \cite{dap11kent} in 2011.  He also did a small amount of teaching at QMUL: for example, each year he gave two undergraduate lectures on the design of questionnaires.

On the research front, Donald really switched his base to QMUL, where he gradually moved himself from the Design of Experiments seminar to the Combinatorics Study Group, where he was a regular speaker.  He would talk about mathematics to anyone who would listen, treating everyone as an equal. A typical reaction would be ``Donald has invented a new concept, given me seventeen examples and a conjecture: now my job is to define the concept and prove the conjecture.''

Just as at East Malling, Donald became involved in the QMUL Luncheon Club. This surprised other QMUL staff, because, although the club is open to all past and present staff of all grades, most non-retired staff are either unaware that they may join it or believe that they do not have time to attend. Through the Luncheon Club he met people from several different departments at QMUL.  In particular, he made contact with the people in charge of the organ in the Great Hall and was given permission to play it.

What surprised people much more was Donald's decision to take up 
rock-climbing.  He joined a club in Mile End Park, and learnt to climb on their climbing wall.  Thereafter, he would astound Byron Morgan with announcements like ``I shall take a small diversion on my way back from Edinburgh so that I can climb up the outside of such-and-such an industrial chimney.'' Then he progressed to abseiling down the outside of large buildings to raise money for charity.  He even started to take up sky-diving, but had to stop all of this activity because of a detached retina. Nevertheless, he never displayed any sign of external frailty, and he remained vigorous and upright.

One strand of Donald's research after 2000 concerned neighbour-balanced designs. He used the idea of a \textit{terrace} that I had defined while at Rothamsted.  This is a row of different numbers arranged in a special way: what matters is which number is next door to which other number.  This made me think of a terrace of houses: hence the name.  Unfortunately, one of his ODA colleagues assumed that I meant designs suitable for use in experiments on terraces on steep hillsides in Nepal.  However, by 2000 Donald seemed happy with the word, and produced a string of papers, some with my PhD student Matt Ollis, giving new constructions for terraces. Talks and papers had attractive titles like \textit{Dancing on Ramsgate sands}.

He also started a new fruitful collaboration on this with Ian Anderson of the University of Glasgow.  Ian and Donald had 
lost touch after St Andrews but
met again at one of the BCCs.  Ian was approached by his statistical colleague Tom
Aitchison to construct a schedule for treatment of groups of patients which balanced consecutive pairs of treatments in a particular way.  Ian wrote a draft, and sent it to Donald, asking if the idea was new. He replied that he thought it was, but suggested expanding it, and the result
was a joint paper \cite{dap02and} on locally balanced change-over designs in \textit{Utilitas Mathematica} in 2002.  This is yet another generalization of the rectangular layout: patients are columns, periods are rows, and the columns are clumped into groups.  The balance on consecutive pairs of treatments is akin to, but not the same as, the neighbour-balance property of a quasi-complete Latin square.

Thus began a study of terraces that lasted a decade and produced over a dozen
joint papers.  The basic idea was to construct terraces in the set of integers modulo~$n$ by fitting together sequences of powers of certain integers in clever ways: they called these ``power sequences''. Donald later suggested the most unlikely, but fruitful, strategy of using the arithmetic modulo~$n$ to construct terraces for the integers modulo $m$ where $m=n+1$ or $n-1$ or $n-2$.  Donald could produce endless supplies of elegant examples, and it was Ian's task to provide the number  theory which formed their basis.

Before this collaboration, Ian had worked a lot on whist tournaments.  Now Donald realised that these formed examples of the nested designs that he also studied.

When Donald learnt that Ian's children were musically inclined, he showed great interest, remembering the family chamber music sessions of his childhood.  On his Edinburgh visits, he would go across to Glasgow to play chamber music with Ian's children.  He passed on to them scores of various trios and quartets.


Another long collaboration was with Peter Cameron.  When the number of treatments is a prime power, such as $4$ or $13$, there are often general constructions that can be made from a single starting array by using a primitive element of the finite field of that order.  Things are more difficult for composite numbers.  In terminology from R.~D.~Carmichael
\cite{car1,car2},
a \textit{primitive lambda-root} modulo an integer~$n$ is an element of maximal order in the multiplicative group of units 
modulo~$n$.  Now the trick was to find several starting arrays that could be expanded by using a primitive lambda-root and then glued together.  Donald was an expert at finding suitable pieces to glue together, while Peter developed more general theory.  Together they produced an ever-growing series of lecture notes \cite{dap14pjc}, too long for a paper, but too short for a book.

As Peter said, in these collaborations it was as if Donald had climbed straight up the wall by instinct, leaving us to put in the pitons and mark the way so that others could follow.

In 2004, Donald decided to arrange a celebration for his 65th birthday. It was held at QMUL in November.  The formal part consisted of music and mathematics, with Donald contributing to both.  This was followed by a buffet supper and drinks, for mathematicians, musicians and family.
As part of the preparation for this event, Donald managed to re-establish contact with some family members whom he had not met for decades.  Interaction with them gave him much pleasure in the following years.

\section{Retirement}

In July 2008, Donald was told that he would have to retire from QMUL.  At that time, the official retirement age was 65, which he had already passed.
He asked to defer retirement, and even expand his teaching, but this was not allowed. Instead, he was granted emeritus professorship in September 2008.  The accompanying letter stated that this was ``the last  such title approved by Professor Adrian Smith in his role as Principal of the college, which gave great pleasure to both AFMS and DAP.''

Donald now had to vacate his office and move to the room shared by emeritus staff.  Colleagues were handed books and reprints, some signed ``D.~A.~Preece'', others signed ``F.~Yates''.  He continued to come into QMUL two or three days per week,  but he began to have health problems and stays in hospital.  He regaled colleagues with the details: we did not mind; we took the view that we were his family and he needed to talk to someone.  Sometimes he threatened to die before the next week, but there was always a twinkle in his eye and so no one believed him.

In October 2009 he laid on another birthday party at QMUL.  A session on ``Reminiscences and music'' in the Octagon was followed by a buffet in the Senior Common Room.
He also marked his 70th birthday with a special lecture on mathematics in literature at the University of Kent, followed by a piano recital.  He repeated the lecture at QMUL in 2010, and turned the material into a 
paper \cite{dap12hum} published in 2012.

In August 2010 Donald sent an email entitled `triple Youden rectangles' to six co-authors spread around the world.  It really shows his ongoing enthusiasm and excitement about the mathematics, as well as his typical comments on health, so I quote it in full.
\begin{quotation}
Mirabile dictu, I seem to have found the two infinite series that
I have wanted for such a long time.   The `truc' for forcing the
construction to work in general (the multiplication by $2^{-1} \bmod k$)
is so simple that I can only squirm in embarrassment that I didn't
spot it before.   Ah, hindsight\ldots!

The whole thing will have to be written up from scratch, probably
with at least one of you helping me with proofs.   In the meantime,
I've merely revised the text that some of you have seen before,
and a \texttt{.dvi} copy of it is attached.   I fear that it makes very
heavy reading, but I can see how to rewrite the thing in a way
that will make for much easier understanding.

To have achieved this after so long is exciting to the point that
I became worried about my blood-pressure!   However, the readings
that I took were `normal', so it looks like I'll survive this
latest development.

Awkwardly, this brings me to the end of the list of mathematical
things that I said I wanted to see solved ``before I die''!   But I'd
better not say my \textit{Nunc Dimittis} until I've written it up properly!
(If I don't survive that long, any of you should feel free to
complete the task.)

Any comments, thoughts, ideas on the work would of course be very
welcome.

My best wishes to you all.

Donald Preece
\end{quotation}

\textit{Mirabile} indeed, it was John Morgan who replied, even though it was nearly a decade since they had last worked together.  They started work on another joint paper \cite{dap14jpm}, this time on multi-layered Youden rectangles, exchanging versions by email, and meeting during John's 2011 visit to the UK. John would write a theorem or two; Donald would respond that these were not quite covering all that he thought should be possible. With each iteration Donald would explain more, until finally John felt that he understood.  Unfortunately, John became bogged down in university administration, as so many of us do,
and the nearly-complete paper was unfinished at the time of Donald's death.  Only a few months earlier he had sent an email including the words ``I'd like it published in my lifetime, but if you ever hear that I'm no longer around, do proceed with it without me.'' 

In June 2011 Donald finally gave up all teaching at the University of Kent.
He composed an organ piece for the mathematics graduation ceremonies held in Canterbury cathedral the following month.  It was called \textit{Academic fanfare and verse}, and was based on the student song beginning 
``Gaudeamus igitur, juvenes dum sumus''.

Donald then gave himself an interesting new project.  He set himself the task of finding out about all of the pipe organs in the East End of London.  He visited them; he played them; he talked to vicars, organists and historians; he consulted archives.  He put all of this knowledge into a book entitled \textit{The Pipe-Organs of London's East End and its People's Palaces}, which he typeset himself in \LaTeX: it was published by QMUL in 2013.
The book \cite{dap13book} is selling well. It had a very nice review recently in the
\textit{British Journal of Organ Studies}, where he was warmly welcomed as a new author in the field.

Of course, one of the East End organs was the Rutt organ in the Great Hall of the People's Palace at QMUL, which had begun its life bringing music to the local people in the first half of the twentieth century.  Donald was afraid that this organ would be destroyed in the refurbishment of the Great Hall.  He liaised with Philip Ogden, the Senior Adviser to Principal Simon Gaskell,  previously Senior Vice-Principal with special responsibility for estate development and also acting Principal for the year 2008--2009. Through them, he persuaded QMUL to restore the organ.  A concert was held in the Great Hall in March 2013 to inaugurate the refurbished organ: Donald was very proud not only to give a speech at this but to hear the first public performance of the piece \textit{Nostalgic interlude} that he had written especially for this organ. The organist, Alan Wilson, told the audience that it featured ``beautiful undulating string sounds alongside solo flutes''.

Donald was one of the invited speakers at the conference 
\textit{Combinatorics, Algebra and More} held at QMUL in July 2013 to celebrate Peter Cameron's retirement.  He was on good form, lively and cheerful.  His talk used the tredoku puzzle from one of the daily newspapers and turned it into a mathematical question.  During another speaker's later talk, which was admittedly rather opaque, I was amused to see that audience members who had given up attempting to follow it had instead started trying to work on Donald's question.

Donald had a chesty cough for most of 2013.  Although he complained about it, those of us who saw him weekly did not think that he was getting any worse.  In November he was one of the invited speakers at the \textit{Old Codger's} combinatorics colloquium in Reading organized by Anthony Hilton.  Some attenders reported that he was looking poorly; Ian Anderson said that Donald said things like ``I do not think that I shall be able to finish this piece of work, so some of you will have to do it.'' However, he was back on form in QMUL in December, talking excitedly about more East End organs that he had discovered since finishing his book, and passing on from Gavin Ross a jokey verse about R.~A.~Fisher written by his non-statistical colleagues at Rothamsted at Christmas 1919.
His last composition was an organ piece written for the QMUL graduation ceremony in December 2013 in the Great Hall.

As usual, Donald spent the Christmas/New Year period at the house in Edinburgh.  While there, he was taken to the Western General Hospital with a viral infection. He started to respond to treatment, but then died of an apparently unrelated brain haemorrhage.

When the news was passed to people in the School of Mathematical Sciences at QMUL, the responses were remarkably similar.  From PhD students to professors emeriti, from statisticians and astronomers to pure and applied mathematicians, people said ``How can this be? He was still doing so much.  He was my friend.''  We are all stunned, but I think that we are also all grateful that he never had the `retirement with nothing to do' that he so dreaded.

In John Morgan's words: ``How healthy and vital Donald always appeared, and how much energy he exuded, despite his ceaseless litany of health complaints.  I will think of his demise as a result of his brain being no longer able to contain his endless mental energy.'' 

\paragraph{Acknowledgements}
I should like to thank the following for information and reminiscences:
Ian Anderson, Chris Brien, Peter Cameron, Bernard Carr, Robin Constable,
Richard Cormack, Steven Gilmour, John Gower, Mariana Iossifova-Kelly,  Wojtek Krzanowski, John Lamb,
Sue Lewis, Byron Morgan, John Morgan, Philip Ogden, Martin Owen,
Robert Preece, Martin Ridout, Gavin Ross and Iwan Williams

\renewcommand{\refname}{Publications of D.~A.~Preece}

\renewcommand{\refname}{Other references}

\end{document}